# SOME CHARACTERIZATIONS OF QUATERNIONIC RECTIFYING CURVES IN THE SEMI-EUCLIDEAN SPACE $\mathbb{E}_2^4$

Tülay SOYFİDAN and Mehmet Ali GÜNGÖR

**Abstract:** The notion of rectifying curve in the Euclidean space is introduced by Chen as a curve whose position vector always lies in its rectifying plane spanned by the tangent and the binormal vector field $t$ and $n_2$ of the curve. In this study, we have obtained some characterizations of semi-real spatial quaternionic rectifying curves in $\mathbb{R}_1^3$. Moreover, by the aid of these characterizations, we have investigated semi real quaternionic rectifying curves in semi-quaternionic space $\mathbb{Q}_v$.

## 1. Introductions

Quaternions were discovered by Irish mathematician Sir William R. Hamilton who invented quaternions in order to extend three dimensional vector algebra for inclusion of multiplication and divisions. Although it has been found that ordinary vector algebra provides a better mathematical apparatus for investigating physical problems, quaternion algebra, nevertheless, provides us with a simple and elegant represantation for describing finite rotations in space. Moreover, quaternions are found applications in robotics, molecular modelling and computer graphics.

A striking feature of quaternions is that the product of two quaternions is non-commutative, meaning that the product of two quaternions depends on which factor is to the left of the multiplication sign and which factor is to the right. At the beginning, the pratical use of quaternions was minimal in comparison with other methods. But today, quaternions have played a significant role recently in several areas of physical science; namely, in differential geometry, in analysis and synthesis of mechanism and machines simulation of particle motion in molecular physics and quaternionic formulation of equation of motion in the theory of relativity. Moreover, quaternions have developed a wide-spread use in computer graphics and robotics research because they can be used to control rotations in three dimensional space.

As a set, the quaternions $\mathbb{Q}$ are coincide with $\mathbb{R}^4$, a $4$-dimensional vector space over the real numbers. Considering this feature of quaternions, Bharathi and Nagaraj introduced the Serret-Frenet formulae for a quaternionic curves in $\mathbb{R}^3$ and $\mathbb{R}^4$, [2]. Then, lots of papers have been published by using these studies. One of them is a study given by Çöken and Tuna. They gave the Serret-Frenet formulas, inclined curves, harmonic curvatures and some characterizations for a quaternionic curve in the semi-Euclidean spaces $\mathbb{E}_1^3$ and $\mathbb{E}_2^4$, [3].

On the other hand, the notion of the rectifying curve is introduced in [1] as a space curve whose position vector always lies in its rectifying plane. Therefore, the position vector with respect to a chosen origin, of a rectifying curve $\alpha$ in $\mathbb{R}^3$, satisfies the equation
$$\alpha(s) = \lambda(s)t(s) + \mu(s)n_2(s)$$
where $\lambda$ and $\mu$ are arbitrary differentiable functions in terms of the arc length parameter $s \in I \subset \mathbb{R}$.

Chen and Dillen [5] found a relationship between the rectifying curves and the centrodes given by the endpoints of the Darboux vector of a space curve and playing an important role in mechanics. In [6], Ilarslan and Nesovic, defined a rectifying curve in $\mathbb{R}^4$ as a curve whose position vector always lies in the orthogonal complement $N_1^\perp$ of its principal normal vector field $N_1$. Thus, the position vector with respect to a chosen origin, of a rectifying curve in $\mathbb{R}^4$, satisfies the equation
$$\beta(s) = \lambda(s)T(s) + \mu(s)N_2(s) + \nu(s)N_3(s)$$
where $\lambda$, $\mu$ and $\nu$ are differentiable functions in the arc length function $s$.

In addition that, the rectifying curves in Minkowski space were studied by Ilarslan and Nesovic, [7,10].

Then, Güngör and Tosun defined the spatial quaternionic rectifying curves in $\mathbb{R}^3$ and gave some characterizations for these curves. Moreover, they investigated quaternionic rectifying curves in $\mathbb{R}^4$, [4].

The main purpose of this paper, firstly, is to establish some characterizations of spatial quaternionic rectifying curves in $\mathbb{R}_1^3$ and, secondly, to give quaternionic rectifying curves in semi-quaternionic space $\mathbb{Q}_v$ by the aid of these characterizations.

## 2. Preliminaries

A semi-real quaternion is defined with $q = q_1 e_1 + q_2 e_2 + q_3 e_3 + q_4$ (or $q = S_q + V_q$ where the symbols $S_q = q_4$ and $V_q = q_1 e_1 + q_2 e_2 + q_3 e_3$ denote scalar and vector part of $q$, respectively) such that

$$\begin{aligned}
i) \quad & e_i \times e_i = -\varepsilon_{e_i}, \quad 1 \leq i \leq 3, \quad \varepsilon_{e_i} = \pm 1 \\
ii) \quad & e_i \times e_j = \varepsilon_{e_i} \varepsilon_{e_j} e_k \quad \text{in } R_1^3 \\
& e_i \times e_j = -\varepsilon_{e_i} \varepsilon_{e_j} e_k \quad \text{in } R_2^4
\end{aligned} \qquad (2.1)$$

where $(ijk)$ is an even permutation of $(123)$, [3].

For every $p, q \in \mathbb{Q}_v$, using these basic products we can expand the product of two semi-real quaternion as

$$p \times q = S_p S_q + g(V_p, V_q) + S_p V_q + S_q V_p + V_p \wedge_L V_q, \qquad (2.2)$$

where we have used the usual inner and cross products in semi-Euclidean space $\mathbb{R}_1^3$, [3].

The conjugate of the semi-real quaternion $q$ is denoted by $\bar{q}$ and defined $\bar{q} = S_q - V_q$. Thus, we define symmetric, non-degenerate valued bilinear form $h$ as follows:

$$\begin{aligned} h(p,q) &= \frac{1}{2}\left[\varepsilon_p \varepsilon_q (p \times q) + \varepsilon_q \varepsilon_p (q \times p)\right] \quad \text{for } \mathbb{R}_1^3 \\ h(p,q) &= \frac{1}{2}\left[-\varepsilon_p \varepsilon_q (p \times q) - \varepsilon_q \varepsilon_p (q \times p)\right] \quad \text{for } \mathbb{R}_2^4 \end{aligned} \qquad (2.3)$$

and it is called semi-real quaternion inner product, [3].

The norm of a semi-real quaternion $q = (q_1, q_2, q_3, q_4) \in Q_v$ is

$$N(q) = \sqrt{\left| q_1^2 + q_2^2 - q_3^2 - q_4^2 \right|}. \qquad (2.4)$$

If $N(q) = 1$, then $q$ is called a semi-real unit quaternion, [8].

$q$ is called a semi-real spatial quaternion whenever $q + \bar{q} = 0$, [3]. Moreover, quaternion product of two semi-real spatial quaternions is

$$p \times q = g\langle p, q \rangle + p \wedge_L q. \qquad (2.5)$$

$q$ is a semi-real temporal quaternion whenever $q - \bar{q} = 0$. Any general $q$ can be written as

$$q = \frac{1}{2}(q + \bar{q}) + \frac{1}{2}(q - \bar{q}).$$

3. **Some Characterizations of Semi-Real Spatial Quaternionic Curves**

The 3−dimensional semi-Euclidean space $\mathbb{R}_1^3$ is identified with the space of spatial quaternions $\{q \in Q_v \mid q + \bar{q} = 0\}$ in an obvious manner. Let $I = [0,1]$ be an interval in the real line $\mathbb{R}$ and

$$\alpha : I \subset R \to Q_v, \qquad s \to \alpha(s) = \sum_{i=1}^{3} \alpha_i(s) e_i, \qquad 1 \leq i \leq 3,$$

be the smooth curve with the arc-length parameter $s$. So, the tangent $\alpha'(s) = t$ has unit length, [3].

**Theorem 3.1.** Let $\{t, n_1, n_2\}$ be the Serret-Frenet frame at the point $\alpha(s)$ of the semi-real spatial quaternionic curve $\alpha$. Then the Serret-Frenet equations are

$$t' = \varepsilon_{n_1} k n_1,$$
$$n_1' = -\varepsilon_t k t + \varepsilon_{n_1} r n_2, \qquad (3.1)$$
$$n_2' = -\varepsilon_{n_2} r n_1,$$

where $t, n_1, n_2, k$ and $r$ denote the unit tangent, the principal normal, the unit binormal, the principal curvature and the torsion of the semi-real spatial quaternionic curve $\alpha$, respectively. Moreover, $h(t,t) = \varepsilon_t$, $h(n_1, n_1) = \varepsilon_{n_1}$, $h(n_2, n_2) = \varepsilon_{n_2}$, [3].

In a similar manner to the reference [4], we can define the semi-real spatial quaternionic rectifying curves as follows. The position vector of the semi-real spatial quaternionic rectifying curve satisfies the following equation

$$\alpha(s) = \lambda(s) t(s) + \mu(s) n_2(s)$$

where $\lambda$ and $\mu$ are arbitrary differentiable functions.

The following theorems provide some simple characterizations of semi-real spatial quaternionic rectifying curves.

**Theorem 3.2.** Let $\alpha : I \subset R \to R_1^3$ be a semi-real spatial quaternionic curve in $\mathbb{R}_1^3$ with $k > 0$ and let $s$ be its pseudo arc length function of curve. Then $\alpha$ is semi-real spatial quaternionic rectifying curve if and only if the following four statements hold;

  *i*) The distance function $\rho = N(\alpha)$ provides $\rho^2 = |\varepsilon_t s^2 + c_1 s + c_2|$, for some constants $c_1$ and $c_2$.

  *ii*) The tangential component of the position vector of the $\alpha$ is given by $h(\alpha, t) = \varepsilon_t s + c$ for some constant $c$, where $h(t,t) = \varepsilon_t$.

  *iii*) The normal component $\alpha^N$ of the position vector of the $\alpha$ has constant length. In addition, the distance function $\rho$ is non-constant.

  *iv*) The torsion $r(s) \neq 0$ and the binormal component of the position vector $\alpha(s)$ is constant i.e., $h(\alpha, n_2)$ is constant.

**Proof.** Without loss of the generality, let us first suppose that $\alpha : I \subset R \to R_1^3$ is parameterized by the pseudo arc length function $s$ and $\alpha$ is a semi-real spatial quaternionic rectifying curve. So, the position vector $\alpha(s)$ of the $\alpha$ satisfies the equation

$$\alpha(s) = \lambda(s) t(s) + \mu(s) n_2(s) \qquad (3.2)$$

where $\lambda(s)$ and $\mu(s)$ are some differentiable functions of the pseudo arc length parameter $s$.

Differentiating the equation (3.2) with respect to the pseudo arc length function $s$ and by applying the Frenet equation (3.1), we obtain that

$$h(\alpha',t) = \lambda'(s) = 1$$
$$h(\alpha',n_1) = \lambda(s)k(s) - \varepsilon_{n_1}\varepsilon_{n_2}\mu(s)r(s) = 0 \qquad (3.3)$$
$$h(\alpha',n_2) = \mu'(s) = 0$$

whereby $h(n_1,n_1) = \varepsilon_{n_1}$ and $h(n_2,n_2) = \varepsilon_{n_2}$.

Moreover, we get

$$\lambda(s) = s + c,$$
$$\mu(s) = a, \qquad (3.4)$$
$$\lambda(s)k(s) = \varepsilon_{n_1}\varepsilon_{n_2}\mu(s)r(s) \neq 0,$$

where $c, a \in R$ and hence $a \neq 0$, $r(s) \neq 0$. From the equation (3.2) we easily obtain that

$$\rho^2 = N^2(\alpha) = |h(\alpha,\alpha)| = |\varepsilon_t \lambda^2 + \varepsilon_{n_2}\mu^2|.$$

By substituting the equation (3.4) into the last equation, we get $\rho^2 = |\varepsilon_t s^2 + c_1 s + c_2|$. This proves statement $(i)$.

Further, from the equations (3.2) and (3.4) we get $h(\alpha,t) = \varepsilon_t \lambda = \varepsilon_t s + c$. This implies statement $(ii)$.

Next, from the equation (3.2) we have the normal component $\alpha^N$ of the position vector of the semi-real spatial quaternionic rectifying curve is given by $\alpha^N = \mu n_2$. Therefore, $N(\alpha^N) = |a| \neq 0$. Thus, statement $(iii)$ is proved.

And lastly, from (3.2) we easily find that $h(\alpha,n_2) = \mu = $ constant and since $r(s) \neq 0$, the statement $(iv)$ is proved.

Conversely, assume that statement $(i)$ or statement $(ii)$ holds. Then, we have $h(\alpha,t) = \varepsilon_t s + c$, for some constant $c$. Differentiating this equation with respect to $s$, we obtain that $\varepsilon_{n_1} k h(\alpha,n_1) = 0$. Hence $h(\alpha,n_1) = 0$ is found. So, $\alpha$ is a semi-real spatial quaternionic rectifying curve.

Next, if we accept the statement $(iii)$ holds, then we have $\alpha(s) = m(s)t(s) + \alpha^N(s)$, $m(s) \in R$. Hence we easily find that

$$c = h(\alpha^N, \alpha^N) = h(\alpha,\alpha) - \varepsilon_t h(\alpha,t)^2.$$

By differentiating the last equation with respect to $s$ gives

$$h(\alpha,t) = \varepsilon_t h(\alpha,t)\left(\varepsilon_t + \varepsilon_{n_1} k\, h(\alpha,\mathbf{n}_1)\right). \tag{3.5}$$

On the other hand, the distance function $\rho = N(\alpha)$ is a non-constant. So, we get $h(\alpha,t) \neq 0$. Therefore, since function $k(s) > 0$ and from (3.5) we find $h(\alpha,\mathbf{n}_1) = 0$, which means that $\alpha$ is a semi-real spatial quaternionic rectifying curve.

And finally if statement (iv) holds, we have $h(\alpha,\mathbf{n}_2) = \mu = $ constant. If we take the derivative of this equation with respect to $s$ and apply Frenet equations (3.1), then we easily obtain $\varepsilon_{n_2} r h(\alpha,\mathbf{n}_1) = 0$. Since $r(s) \neq 0$, we find that $h(\alpha,\mathbf{n}_1) = 0$. So, it means that the curve $\alpha$ is a semi-real spatial quaternionic rectifying curve.

In the next theorem, we prove that the ratio of the torsion and curvature of a unit speed semi-real spatial quaternionic rectifying curve is a non-constant linear function of the pseudo arc length parameter $s$.

**Theorem 3.3.** Let $\alpha : I \subset R \to R_1^3$ be a semi-real spatial quaternionic curve in $\mathbb{R}_1^3$ with a spacelike or a timelike rectifying plane and with the curvature $k(s) > 0$. The curve $\alpha$ is a semi-real spatial quaternionic curve if and only if there holds $\dfrac{r(s)}{k(s)} = c_1 s + c_2$, where $c_1, c_2 \in R$, $c_1 \neq 0$.

**Proof.** Without loss of generality, suppose that $\alpha : I \subset \mathbb{R} \to \mathbb{R}_1^3$ is parameterized by the pseudo arc length function $s$ and $\alpha$ is a semi-real spatial quaternionic rectifying curve. By the use of the equation (3.4), we have $\dfrac{r(s)}{k(s)} = \varepsilon_{n_1}\varepsilon_{n_2} \dfrac{\lambda(s)}{\mu(s)} = \varepsilon_{n_1}\varepsilon_{n_2} \dfrac{s+c}{a}$ for some constants $c_1$ and $a$. Consequently, $\dfrac{r(s)}{k(s)} = c_1 s + c_2$, whereby $c_1 \neq 0$ and $c_1, c_2 \in R$.

Conversely, let us suppose that $\dfrac{r(s)}{k(s)} = c_1 s + c_2$, $c_1, c_2 \in R$ and $c_1 \neq 0$. Then, especially we take $c_1 = \dfrac{\varepsilon_{n_1}\varepsilon_{n_2}}{a}$, $c_2 = \dfrac{\varepsilon_{n_1}\varepsilon_{n_2} c}{a}$, where $c \in R$, $a > 0$. So, $\dfrac{r(s)}{k(s)} = \varepsilon_{n_1}\varepsilon_{n_2}\dfrac{s+c}{a}$. Thus, by using the Frenet equations (3.1), we obtain that

$$\frac{d}{ds}\left[\alpha(s) - \lambda(s)\mathbf{t}(s) - \mu(s)\mathbf{n}_2(s)\right] = 0,$$

which means that the curve $\alpha$ is semi-real spatial quaternionic rectifying curve.

**Theorem 3.4.** Let $\alpha = \alpha(s)$ be a unit speed semi-real spatial quaternionic curve in $\mathbb{R}_1^3$. Then the following statements hold:

*i)*     $\alpha$ is a semi-real spatial quaternionic rectifying curve with a spacelike rectifying plane if and only if, up to parametrization, $\alpha$ is given by

$$\alpha(t) = y(t)\frac{a}{\cos t}, \quad a > 0,$$

where $y(t)$ is a unit speed spatial spacelike quaternionic curve lying in the pseudosphere $S_1^2(1)$.

*ii)*     $\alpha$ is a semi-real spatial spacelike (timelike) quaternionic rectifying curve with a timelike rectifying plane and a spacelike (timelike) quaternionic position vector if and only if, up to parametrization, $\alpha$ is given by

$$\alpha(t) = y(t)\frac{a}{\sinh t}, \quad a > 0,$$

where $y(t)$ is a unit speed spatial timelike (spacelike) quaternionic curve lying in the pseudosphere $S_1^2(1)$ (pseudohyperbolic space $H_0^2(1)$).

*iii)*     $\alpha$ is a semi-real spatial spacelike (timelike) quaternionic rectifying curve with a timelike rectifying plane and a timelike (spacelike) quaternionic position vector if and only if, up to parametrization, $\alpha$ is given by

$$\alpha(t) = y(t)\frac{a}{\cosh t}, \quad a > 0,$$

where $y(t)$ is a unit speed spatial spacelike (timelike) quaternionic curve lying in the pseudohyperbolic space $H_0^2(1)$ (pseudosphere $S_1^2(1)$).

**Proof:** *i)* Let $\alpha(s)$ be a unit speed semi-real quaternionic rectifying curve with spacelike rectifying plane in $\mathbb{R}_1^3$. Since the position vector lies in the spacelike rectifying plane, we have $h(\alpha, \alpha) > 0$, $h(t,t) = \varepsilon_t = 1$ and $h(n_2, n_2) = \varepsilon_{n_2} = 1$. By the proof of Theorem 3.2, we have $\rho^2 = N(\alpha)^2 = (s+c)^2 + a^2$, $s, a \in R$, $a > 0$. If we apply a translation with respect to $s$, then we have $\rho^2 = s^2 + a^2$.

Next, we define a curve $y$ lying in the pseudosphere $S_1^2(1)$ by

$$y(s) = \frac{\alpha(s)}{\rho(s)}. \tag{3.6}$$

Then, we get

$$\alpha(s) = y(s)\sqrt{s^2 + a^2}. \tag{3.7}$$

By differentiating equation (3.7) with respect to $s$, we obtain that

$$t(s) = y(s)\frac{s}{\sqrt{s^2+a^2}} + y'(s)\sqrt{s^2+a^2}. \tag{3.8}$$

Since $h(y,y)=1$, we find that $h(y,y')=0$. From equation (3.8) we get

$$1 = h(t,t) = h(y',y')(s^2+a^2) + \frac{s^2}{s^2+a^2},$$

and so

$$h(y',y') = \frac{a^2}{\left(s^2+a^2\right)^2}, \tag{3.9}$$

which means that $y$ is a spacelike semi-real quaternionic curve.

From (3.9) we obtain that $N(y'(s)) = \frac{a}{s^2+a^2}$. Let $t = \int_0^s N(y'(u))du$ be the pseudo arc length parameter of the semi-real quaternionic curve $y$. Then we have

$$t = \int_0^s \frac{a}{u^2+a^2}du,$$

and therefore $s = a\tan t$. By substituting this into equation the (3.7) we obtain the parametrization.

Conversely, assume that $\alpha$ is a semi-real quaternionic curve defined by (i), where $y(t)$ is a unit speed semi-real spacelike quaternionic curve lying in the pseudosphere $S_1^2(1)$. By differentiating the equation (i) with respect to $t$, we get

$$\alpha'(t) = \frac{a}{\cos^2 t}\left(y(t)\sin t + y'(t)\cos t\right).$$

Moreover, we have $h(y',y')=1$, $h(y,y)=1$ and $h(y,y')=0$. So, we obtain that

$$h(\alpha,\alpha') = \frac{a^2 \sin t}{\cos^3 t}, \quad h(\alpha',\alpha') = \frac{a^2}{\cos^4 t}, \tag{3.10}$$

and consequently $N(\alpha'(t)) = \frac{a}{\cos^2 t}$. Let us put $\alpha(t) = m(t)\alpha'(t) + \alpha^N$, where $m(t) \in R$ and $\alpha^N$ is a normal component of the position vector $\alpha$. Then we easily find that $m = \frac{h(\alpha,\alpha')}{h(\alpha',\alpha')}$, and therefore

$$h(\alpha^N,\alpha^N) = h(\alpha,\alpha) - \frac{h(\alpha,\alpha')^2}{h(\alpha',\alpha')}.$$

Since $h(\alpha,\alpha) = \dfrac{a^2}{\cos^2 t}$ and by using (3.10), the last equation converts into $h(\alpha^N, \alpha^N) = a^2 =$ constant. So, we find that $N(\alpha^N(t)) =$ constant and since $\rho = N(\alpha) = \dfrac{a}{\cos t} \neq$ constant, Theorem 3.2 implies that $\alpha$ is a semi-real quaternionic rectifying curve.

*ii)* Let $\alpha(s)$ be a unit speed semi-real spacelike quaternionic rectifying curve with a timelike rectifying plane and a spacelike position vector. Hence we have $h(\alpha,\alpha) > 0$, $h(t,t) = \varepsilon_t = 1$ and $h(n_2,n_2) = \varepsilon_{n_2} = -1$. By the proof of Theorem 3.2, we obtain $\rho^2 = N(\alpha)^2 = (s+c)^2 - a^2$, $s, a \in R$, $a > 0$. If we apply a translation with respect to $s$, then we have $\rho^2 = s^2 - a^2$, $|s| > a$. Further, we define a curve $y(s)$ lying in the pseudosphere $S_1^2(1)$ by

$$y(s) = \frac{\alpha(s)}{\rho(s)}. \tag{3.11}$$

It follows that

$$\alpha(s) = y(s)\sqrt{s^2 - a^2}, \tag{3.12}$$

and if we take the derivative of the previous equation with respect to $s$, then we easily find that

$$t(s) = y(s)\frac{s}{\sqrt{s^2 - a^2}} + y'(s)\sqrt{s^2 - a^2}. \tag{3.13}$$

Since $h(y,y) = 1$, it follows that $h(y,y') = 0$.
Consequently,

$$h(t,t) = h(y',y')(s^2 - a^2) + \frac{s^2}{s^2 - a^2} = 1.$$

From the last equation if we do necessary arrangement, we obtain that

$$h(y',y') = \frac{a^2}{(s^2 - a^2)^2}, \tag{3.14}$$

which means that $y$ is a semi-real timelike quaternionic curve.

By considering equation (3.14), we easily find that $N(y'(s)) = \dfrac{a}{s^2 - a^2}$, $a > 0$, $|s| > a$. Let $t = \displaystyle\int_0^s N(y'(u))\,du$ be the pseudo arc length parameter of the semi-real quaternionic curve $y$. So, we have

$$t = \int_0^s \frac{a}{u^2 - a^2} du,$$

and thus $s = -a \coth t$. By substituting this into (3.12), we obtain the parametrization.

Conversely, let us assume that $\alpha$ is a semi-real quaternionic curve defined by statement (*ii*), where $y(t)$ is a unit speed semi-real timelike quaternionic curve lying in the pseudosphere $S_1^2(1)$. By differentiating the statement (*ii*) with respect to $t$, we find that

$$\alpha'(t) = \frac{a}{\sinh^2 t}\left(y'(t)\sinh t - y(t)\cosh t\right). \tag{3.15}$$

By assumption we have $h(y', y') = -1$, $h(y, y) = 1$ and therefore $h(y, y') = 0$. Then, the equation (3.15) implies that

$$h(\alpha, \alpha') = -\frac{a^2 \cosh t}{\sinh^3 t}, \quad h(\alpha', \alpha') = \frac{a^2}{\sinh^4 t}, \tag{3.16}$$

and therefore $N(\alpha'(t)) = \frac{a}{\cos^2 t}$.

Let us put $\alpha(t) = m(t)\alpha'(t) + \alpha^N$, where $m(t) \in R$ and $\alpha^N$ is a normal component of the position vector $\alpha$. Then we easily find that $m = \frac{h(\alpha, \alpha')}{h(\alpha', \alpha')}$, and therefore

$$h(\alpha^N, \alpha^N) = h(\alpha, \alpha) - \frac{h(\alpha, \alpha')^2}{h(\alpha', \alpha')}. \tag{3.17}$$

Since $h(\alpha, \alpha) = \frac{a^2}{\sinh^2 t}$ and by using the equation (3.16), the equation (3.17) becomes $h(\alpha^N, \alpha^N) = -a^2 = $ constant. Hence $N(\alpha^N(t)) = $ constant and since $\rho = N(\alpha) = \frac{a}{\sinh t} \neq $ constant, Theorem 3.2 implies that $\alpha$ is a semi-real quaternionic rectifying curve.

The proof in the case when $\alpha$ is a semi-real timelike quaternionic rectifying curve with a timelike rectifying plane and a timelike quaternionic position vector is analogous.

*iii*) The proof is analogous to the proofs of the statements *i)* and *ii)*.

## 4 Some Characterizations of Semi-Real Quaternionic Rectifying Curves

In this section, we firstly mention semi-real quaternionic curve in $\mathbb{R}_2^4$. And then we characterize the semi-real quaternonic rectifying curve in $\mathbb{R}_2^4$ in terms of their curvatures.

As in the Section 3, the four dimensional semi-Euclidean space $\mathbb{R}_2^4$ is identified with the space of unit quaternions which is denoted by $\mathbb{Q}_v$. Suppose that

$$\beta: I \subset R \to Q_v, \qquad s \to \beta(s) = \sum_{i=1}^{4} \alpha_i(s) e_i, \quad (1 \le i \le 4), \qquad e_4 = 1$$

is a smooth curve defined over the interval $I = [0,1]$. Let the parameter $s$ be chosen such that the tangent $T = \beta'(s)$ has unit magnitude.

**Theorem 4.1.** Let $\{T, N_1, N_2, N_3\}$ be the Serret-Frenet frame in the point $\beta(s)$ of the semi-real quaternionic curve $\beta$. Then the Serret-Frenet equations are

$$\begin{aligned}
T' &= \varepsilon_{N_1} \kappa N_1, \\
N_1' &= -\varepsilon_t \varepsilon_{N_1} \kappa T + \varepsilon_{n_1} k N_2, \\
N_2' &= -\varepsilon_t k N_1 + \varepsilon_{n_1} \left( r - \varepsilon_t \varepsilon_T \varepsilon_{N_1} \kappa \right) N_3, \\
N_3' &= -\varepsilon_{n_2} \left( r - \varepsilon_t \varepsilon_T \varepsilon_{N_1} \kappa \right) N_2
\end{aligned} \qquad (4.1)$$

where

$$\begin{aligned}
\kappa &= \varepsilon_{N_1} \|T'\|, \ N_1 = \varepsilon_T (t \times T), \ N_2 = \varepsilon_T (n_1 \times T), \ N_3 = \varepsilon_T (n_2 \times T) \\
h(T,T) &= \varepsilon_T, \ h(N_1, N_1) = \varepsilon_{N_1}, \ h(N_2, N_2) = \varepsilon_{n_1} \varepsilon_T, \ h(N_3, N_3) = \varepsilon_{n_2} \varepsilon_T.
\end{aligned} \qquad (4.2)$$

It is obtained the Serret-Frenet formulae and the apparatus for the curve $\beta = \beta(s)$ by making use of the Serret-Frenet formulae for a curve $\alpha = \alpha(s)$ in $\mathbb{R}_1^3$, [3].

Moreover, there are relationships between curvatures of the curves $\beta$ and $\alpha$. These relations can be explained that the torsion of $\beta$ is the principal curvature of the curve $\alpha$. Also, the bitorsion of $\beta$ is $(r - \varepsilon_t \varepsilon_T \varepsilon_{N_1} \kappa)$, where $r$ is the torsion of $\alpha$ and $\kappa$ is the principal curvature of $\beta$, [3].

Now, we characterize the semi-real quaternionic rectifying curves with respect to their curvatures. Let $\beta = \beta(s)$ be unit speed semi-real quaternionic curve in $\mathbb{Q}_v$, with non-zero curvatures $\kappa(s)$, $k(s)$ and $(r - \varepsilon_t \varepsilon_T \varepsilon_{N_1} \kappa)$. Then, the position vector $\beta(s)$ of the semi-real quaternionic curve $\beta$ satisfies the equation

$$\beta(s) = \lambda(s) T(s) + \mu(s) N_2(s) + \nu(s) N_3(s)$$

for some differentiable functions $\lambda(s)$, $\mu(s)$ and $\nu(s)$. Therefore, we give the following theorems of semi-real quaternionic rectifying curve $\beta$.

**Theorem 4.2.** Let us suppose that $\beta = \beta(s)$ is a unit speed semi-real quaternionic curve in $\mathbb{Q}_v$, with non-zero curvatures $\kappa(s)$, $k(s)$ and $(r - \varepsilon_t \varepsilon_T \varepsilon_{N_1} \kappa)$. Then $\beta$ is congruent to a semi-real quaternionic rectifying curve if and only if

$$\varepsilon_t \varepsilon_{n_1} \varepsilon_{N_1} \frac{\kappa(s)(r - \varepsilon_t \varepsilon_T \varepsilon_{N_1} \kappa)(s+c)}{k(s)} + \varepsilon_t \varepsilon_{n_2} \varepsilon_{N_1} \left[ \frac{\kappa(s)k(s) + (s+c)[\kappa'(s)k(s) - \kappa(s)k'(s)]}{k^2(s)(r - \varepsilon_t \varepsilon_T \varepsilon_{N_1} \kappa)} \right]' = 0, c \in R. \quad (4.3)$$

**Proof:** Let $\beta = \beta(s)$ be a unit speed semi-real quternionic rectifying curve and $\kappa(s), k(s)$ and $(r - \varepsilon_t \varepsilon_T \varepsilon_{N_1} \kappa)$ be non zero curvatures of $\beta$. By definition, the position vector of the curve $\beta$ satisfies

$$\beta(s) = \lambda(s)\mathbf{T}(s) + \mu(s)\mathbf{N}_2(s) + \nu(s)\mathbf{N}_3(s) \quad (4.4)$$

for some differentiable functions $\lambda(s)$, $\mu(s)$ and $\nu(s)$.

By differentiating the equation (4.4) with respect to the pseudo arc length function $s$ and using the Frenet equations (4.1), we find that

$$\mathbf{T} = \lambda'\mathbf{T} + (\varepsilon_{N_1}\lambda\kappa - \varepsilon_t\mu k)\mathbf{N}_1 + \left[\mu' - \varepsilon_{n_2}\nu(r - \varepsilon_t \varepsilon_T \varepsilon_{N_1} \kappa)\right]\mathbf{N}_2 + \left[\varepsilon_{n_1}\mu(r - \varepsilon_t \varepsilon_T \varepsilon_{N_1} \kappa) + \nu'\right]\mathbf{N}_3.$$

So, we have

$$\begin{aligned} \lambda' &= 1 \\ \varepsilon_{N_1}\lambda\kappa - \varepsilon_t\mu k &= 0 \\ \mu' - \varepsilon_{n_2}\nu(r - \varepsilon_t \varepsilon_T \varepsilon_{N_1} \kappa) &= 0 \\ \varepsilon_{n_1}\mu(r - \varepsilon_t \varepsilon_T \varepsilon_{N_1} \kappa) + \nu' &= 0 \end{aligned} \quad (4.5)$$

and therefore

$$\begin{aligned} \lambda(s) &= s + c, \\ \mu(s) &= \varepsilon_t \varepsilon_{N_1} \frac{\kappa(s)(s+c)}{k(s)}, \\ \nu(s) &= \varepsilon_t \varepsilon_{n_2} \varepsilon_{N_1} \frac{\kappa(s)k(s) + (s+c)[\kappa'(s)k(s) - \kappa(s)k'(s)]}{k^2(s)(r - \varepsilon_t \varepsilon_T \varepsilon_{N_1} \kappa)(s)}, \end{aligned} \quad (4.6)$$

where $c \in R$.

In this way we express that the functions $\lambda(s)$, $\mu(s)$ and $\nu(s)$ in terms of the curvature functions $\kappa(s)$, $k(s)$ and $(r(s) - \varepsilon_t \varepsilon_T \varepsilon_{N_1} \kappa(s))$ of the semi-real quaternionic rectifying curve $\beta$. Therefore, by using the last equation (4.5) and relation equation (4.6), we find that curvatures $\kappa(s)$, $k(s)$ and $(r(s) - \varepsilon_t \varepsilon_T \varepsilon_{N_1} \kappa(s))$ satisfy the equation (4.3).

Conversely, assume that the curvatures $\kappa(s)$, $k(s)$ and $\left(r(s)-\varepsilon_t\varepsilon_T\varepsilon_{N_1}\kappa(s)\right)$ of an arbitrary unit speed semi-real quaternionic rectifying curve $\beta$ in $\mathbb{Q}_v$. Thus, by considering the Frenet equations (4.1) and (4.3), we find that

$$\frac{d}{ds}\left[\begin{array}{l}\beta(s)-(s+c)T(s)-\varepsilon_t\varepsilon_{N_1}\dfrac{\kappa(s)(s+c)}{k(s)}N_2(s)\\ -\varepsilon_t\varepsilon_{n_2}\varepsilon_{N_1}\dfrac{\kappa(s)k(s)+(s+c)\left[\kappa'(s)k(s)-\kappa(s)k'(s)\right]}{k^2(s)\left(r-\varepsilon_t\varepsilon_T\varepsilon_{N_1}\kappa\right)(s)}N_3(s)\end{array}\right]=0$$

which means that the quaternionic curve $\beta$ is congruent to a semi-real quaternionic rectifying curve.

Now, we assume that all the curvature functions $\kappa(s)$, $k(s)$ and $\left(r(s)-\varepsilon_t\varepsilon_T\varepsilon_{N_1}\kappa(s)\right)$ of quaternionic rectifying curve $\beta$ in $\mathbb{Q}_v$ are constant and different from zero. Then, we can express the following corollary from the Theorem 4.2.

**Corollary 4.1.** There are no quaternionic semi-real rectifying curves lying fully in $\mathbb{Q}_v$, with non-zero constant curvatures $\kappa(s)$, $k(s)$ and $\left(r(s)-\varepsilon_t\varepsilon_T\varepsilon_{N_1}\kappa(s)\right)$.

Let $\beta=\beta(s)$ be a unit speed semi-real quternionic curve in $\mathbb{Q}_v$, with non-zero curvatures $\kappa(s)$, $k(s)$ and $\left(r-\varepsilon_t\varepsilon_T\varepsilon_{N_1}\kappa\right)(s)$. If any two of the curvature functions are constant, we may consider the following theorem.

**Theorem 4.3** Let $\beta=\beta(s)$ be a unit speed semi-real quternionic curve in $\mathbb{Q}_v$, with non-zero curvatures $\kappa(s)$, $k(s)$ and $\left(r-\varepsilon_t\varepsilon_T\varepsilon_{N_1}\kappa\right)(s)$. Then $\beta$ is congruent to a semi-real quaternionic rectifying curve if

i) $\kappa(s)=\text{constant}>0$, $k(s)=\text{constant}\neq 0$ and

$$\left(r-\varepsilon_t\varepsilon_T\varepsilon_{N_1}\kappa\right)=\frac{1}{\sqrt{\left|-\varepsilon_{n_1}\varepsilon_{n_2}s^2-2\varepsilon_{n_1}\varepsilon_{n_2}cs-2c_1\right|}}$$

(namely, $r(s)=\text{non constant}$), $c,c_1\in R$;

ii) $k(s)=\text{constant}\neq 0$ and $r'(s)=\varepsilon_t\varepsilon_T\varepsilon_{N_1}\kappa'(s)$ (namely, $\left(r-\varepsilon_t\varepsilon_T\varepsilon_{N_1}\kappa\right)(s)=\text{constant}$). Then,

$$\kappa(s)=\frac{c_1 e^{i\sqrt{\varepsilon_{n_1}\varepsilon_{n_2}}(r-\varepsilon_t\varepsilon_T\varepsilon_{N_1}\kappa)(s)}+c_2 e^{-i\sqrt{\varepsilon_{n_1}\varepsilon_{n_2}}(r-\varepsilon_t\varepsilon_T\varepsilon_{N_1}\kappa)(s)}}{s+c},\quad c,c_1,c_2\in R.$$

iii) $\kappa(s)=\text{constant}>0$, $r(s)=\text{constant}$ (namely, $\left(r-\varepsilon_t\varepsilon_T\varepsilon_{N_1}\kappa\right)(s)=\text{constant}$) and

$$k(s)=\frac{s+c}{c_1 e^{i\sqrt{\varepsilon_{n_1}\varepsilon_{n_2}}(r-\varepsilon_t\varepsilon_T\varepsilon_{N_1}\kappa)(s)}+c_2 e^{-i\sqrt{\varepsilon_{n_1}\varepsilon_{n_2}}(r-\varepsilon_t\varepsilon_T\varepsilon_{N_1}\kappa)(s)}},\quad c,c_1,c_2\in R.$$

**Proof:**

*i)* Suppose that $\kappa(s) = \text{constant} > 0$, $k(s) = \text{constant} \neq 0$ and $r(s)$ is non-constant function. By using the equation (4.5), we find differential equation

$$r'(s) - \varepsilon_{n_1}\varepsilon_{n_2}\left(r - \varepsilon_t\varepsilon_T\varepsilon_{N_1}\kappa\right)^3 (s+c) = 0, \quad c \in R.$$

The solution of the this differential equation is given by

$$\left(r - \varepsilon_t\varepsilon_T\varepsilon_{N_1}\kappa\right) = \frac{1}{\sqrt{\left|-\varepsilon_{n_1}\varepsilon_{n_2} s^2 - 2\varepsilon_{n_1}\varepsilon_{n_2} cs - 2c_1\right|}}.$$

*ii)* Now, assume that $k(s) = \text{constant} \neq 0$, $r'(s) = \varepsilon_t\varepsilon_T\varepsilon_{N_1}\kappa'(s)$ (namely, $(r - \varepsilon_t\varepsilon_T\varepsilon_{N_1}\kappa)(s) = \text{constant}$) non constant functions. Then equation (4.5) implies differential equation

$$\varepsilon_{n_1}\kappa(s)(s+c)\left[r(s) - \varepsilon_t\varepsilon_T\varepsilon_{N_1}\kappa(s)\right]^2 + \varepsilon_{n_2}\left[\kappa(s)(s+c)\right]'' = 0,\ r(s) - \varepsilon_t\varepsilon_T\varepsilon_{N_1}\kappa(s) = \text{constant} \neq 0,\ c \in R,$$

whose solution has the form

$$\kappa(s) = \frac{c_1 e^{i\sqrt{\varepsilon_{n_1}\varepsilon_{n_2}}(r - \varepsilon_t\varepsilon_T\varepsilon_{N_1}\kappa)(s)} + c_2 e^{-i\sqrt{\varepsilon_{n_1}\varepsilon_{n_2}}(r - \varepsilon_t\varepsilon_T\varepsilon_{N_1}\kappa)(s)}}{s+c}, \quad c, c_1, c_2 \in R.$$

*iii)* If $\kappa(s) = \text{constant} > 0$, $r(s) = \text{constant}$ (namely, $r(s) - \varepsilon_t\varepsilon_T\varepsilon_{N_1}\kappa(s) = \text{constant}$) and $k(s)$ is non-constant function, by using the equation (4.3) we get the following differential equation

$$\varepsilon_{n_1} \frac{\kappa(s)(s+c)\left[r(s) - \varepsilon_t\varepsilon_T\varepsilon_{N_1}\kappa(s)\right]^2}{k(s)} + \left[\frac{\kappa(s)(s+c)}{k(s)}\right]'' = 0,\ r(s) - \varepsilon_t\varepsilon_T\varepsilon_{N_1}\kappa(s) = \text{constant} \neq 0,\ c \in R.$$

The solution of the above differential equation is given by

$$k(s) = \frac{s+c}{c_1 e^{i\sqrt{\varepsilon_{n_1}\varepsilon_{n_2}}(r - \varepsilon_t\varepsilon_T\varepsilon_{N_1}\kappa)(s)} + c_2 e^{-i\sqrt{\varepsilon_{n_1}\varepsilon_{n_2}}(r - \varepsilon_t\varepsilon_T\varepsilon_{N_1}\kappa)(s)}}, \quad c, c_1, c_2 \in R.$$

The following theorem provides some simple characterizations of semi-real quaternionic rectifying curves in $\mathbb{Q}_v$.

**Theorem 4.4.** Let $\beta = \beta(s)$ be a unit speed semi-real quaternionic rectifying curve in $\mathbb{Q}_v$, with non-zero curvatures $\kappa(s)$, $k(s)$ and $(r(s) - \varepsilon_t \varepsilon_T \varepsilon_{N_1} \kappa(s))$ if and only if one the following four statements hold;

  *i)* The distance function $\rho = N(\beta)$ satisfies $\rho^2 = |\varepsilon_T s^2 + c_1 s + c_2|$ for some constants $c_1$ and $c_2$.

  *ii)* The tangential component of the position vector of the semi-real spatial quaternionic curve is given by $h(\beta, T) = \varepsilon_T s + c$ for some constant $c$.

  *iii)* The normal component $\beta^N$ of the position vector of the curve has constant length and the distance function $\rho$ is nonconstant.

  *iv)* The first binormal component and the second binormal component of the position vector of the semi-real quaternionic curve are respectively given by

$$h(\beta(s), N_2(s)) = \varepsilon_t \varepsilon_{n_1} \varepsilon_T \varepsilon_{N_1} \frac{\kappa(s)(s+c)}{k(s)}$$

$$h(\beta(s), N_3(s)) = \varepsilon_t \varepsilon_T \varepsilon_{N_1} \frac{\kappa(s)k(s) + (s+c)[\kappa'(s)k(s) - \kappa(s)k'(s)]}{k^2(s)(r - \varepsilon_t \varepsilon_T \varepsilon_{N_1} \kappa)(s)}, \quad c \in R. \tag{4.7}$$

**Proof:**

  *i)* Let us assume that $\beta : I \to \mathbb{Q}_v$ is parameterized by the pseudo arc length function $s$ and $\beta$ is a semi-real quaternionic rectifying curve in $\mathbb{Q}_v$ with non-zero curvatures $\kappa(s)$, $k(s)$ and $(r - \varepsilon_t \varepsilon_T \varepsilon_{N_1} \kappa)$.

The position vector $\beta(s)$ of the $\beta$ satisfies the equation (4.4), where the differentiable functions $\lambda(s)$, $\mu(s)$ and $\nu(s)$ satisfy the relation (4.5). If we make necessary arrangement in the equation (4.5), we obtain that

$$(r - \varepsilon_t \varepsilon_T \varepsilon_{N_1} \kappa)(s)\left[\varepsilon_{n_2} \mu(s)\mu'(s) + \varepsilon_{n_2} \nu(s)\nu'(s)\right] = 0$$

Since $(r - \varepsilon_t \varepsilon_T \varepsilon_{N_1} \kappa) \neq 0$ by assumption, $\varepsilon_{n_1} \mu(s)\mu'(s) + \varepsilon_{n_2} \nu(s)\nu'(s) = 0$ is found.

So,
$$\varepsilon_{n_1} \mu^2(s) + \varepsilon_{n_2} \nu^2(s) = a^2, \tag{4.8}$$

for some constant $a > 0$. From the equation (4.4) we have $h(\beta(s), \beta(s)) = \varepsilon_T \lambda^2(s) + \varepsilon_{n_1} \varepsilon_T \mu^2(s) + \varepsilon_{n_2} \varepsilon_T \nu^2(s)$, which together with (4.6) and (4.8) gives $h(\beta(s), \beta(s)) = \varepsilon_T \left[(s+c)^2 + a^2\right]$. Therefore, we obtain $\rho^2 = |\varepsilon_T s^2 + c_1 s + c_2|$, for some constants $c_1$ and $c_2$. So, the statement (*i*) is provided.

  *ii)* From the equation (4.4) we obtain that $h(\beta(s), T(s)) = \varepsilon_T \lambda(s)$, which together with (4.6) give $h(\beta(s), T(s)) = \varepsilon_T s + c$, $c \in R$.

  *iii)* From the equation (4.4) it is clear that the normal component $\beta^N$ of the position vector of the quaternionic rectifying curve implies $\beta^N(s) = \mu(s) N_2(s) + \nu(s) N_3(s)$ and therefore $h(\beta^N(s),\beta^N(s)) = \varepsilon_{n_1}\varepsilon_T \mu^2(s) + \varepsilon_{n_2}\varepsilon_T \nu^2(s)$. Therefore, by using (4.8), we find $N(\beta^N(s)) = a$, for some constant $a > 0$. By statement (i), $\rho(s)$ is non-constant function. Thus, statement (iii) is proved.

*iv)* Taking the scalar product of the equation (4.6) with $N_2$ and $N_3$, respectively. Thus, by using the equation (4.6), we find that curvatures $\kappa(s)$, $k(s)$ and $(r(s) - \varepsilon_t \varepsilon_T \varepsilon_{N_1} \kappa(s))$ satisfy the equation (4.7).

Conversely, assume that statement (*i*) (or statement (*ii*)) holds. Then, we have $h(\beta(s),\beta(s)) = |\varepsilon_T s^2 + c_1 s + c_2|$, for some constants $c_1$ and $c_2$. By differentiating this equation two times (one times) with respect to $s$ and using the Frenet equations (4.2), we obtain $\varepsilon_{N_1}\kappa(s)h(\beta(s),N_1(s)) = 0$. Since $\kappa(s) > 0$ by assumption, $h(\beta(s),N_1(s)) = 0$ is found. Hence, $\beta$ is a semi-real quaternionic rectifying curve.

If statement (*iii*) holds, let us put $\beta(s) = \lambda(s)T(s) + \beta^N(s)$, where $\lambda(s)$ is arbitrary differentiable function. Then

$$h(\beta^N(s),\beta^N(s)) = h(\beta(s),\beta(s)) - 2\lambda(s)h(\beta(s),T(s)) + \lambda^2(s)h(T(s),T(s)).$$

Since $h(\beta(s),T(s)) = \varepsilon_T \lambda(s)$ and $h(T(s),T(s)) = \varepsilon_T$, we get

$$h(\beta^N(s),\beta^N(s)) = h(\beta(s),\beta(s)) - \varepsilon_T h(\beta(s),T(s))^2, \tag{4.9}$$

where $h(\beta^N(s),\beta^N(s))$ is a constant from statement (iii) and $h(\beta(s),\beta(s)) = \rho^2(s)$ is a non constant from statement (i).

Differentiating the equation (4.9) with respect to $s$ and using the Serret-Frenet equations (4.4), we find $\varepsilon_T \varepsilon_{N_1} \kappa(s)h(\beta(s),N_1(s)) = 0$. Since $\kappa(s) > 0$ by assumption, $h(\beta(s),N_1(s)) = 0$ is found. Hence, $\beta$ is a semi-real quaternionic rectifying curve.

Finally, if statement (*iv*) holds, by taking the derivative of the first equation in (4.7) with respect to $s$ and using the Serret-Frenet equations (4.2), we obtain

$$-\varepsilon_t k(s)h(\beta(s),N_1(s)) + \varepsilon_{n_1}(r - \varepsilon_t \varepsilon_T \varepsilon_{N_1}\kappa)(s)h(\beta(s),N_3(s)) = \varepsilon_t \varepsilon_{n_1}\varepsilon_T \varepsilon_{N_1}\left[\frac{K(s)(s+c)}{k(s)}\right]'.$$

By using the second equation in (4.7), the last equation becomes $-\varepsilon_t k(s)h(\beta(s),N_1(s)) = 0$. Since $k(s)$ is non zero curvature by assumption, $h(\beta(s),N_1(s)) = 0$ is found. Hence, $\beta$ is a semi-real quaternionic rectifying curve. This proves the theorem.

*Department of Mathematics, Sakarya Unıversity, Sakarya, Turkey.*
*E-mail: tulay.soyfidan@gmail.com*

*Department of Mathematics, Sakarya Unıversity, Sakarya, Turkey.*
*E-mail: agungor@sakarya.edu.tr*